\documentclass{amsart}
\usepackage{amsfonts,amssymb,latexsym,amsthm,amsmath}
\newtheorem{theorem}{Theorem}

\newtheorem{definition}[theorem]{Definition}

\newtheorem*{KM}{Theorem KM}
\newtheorem*{lab}{Theorem L}
\newtheorem*{PH}{Theorem H}

\newtheorem*{P1}{Proposition 1}
\newtheorem*{L1}{Lemma 1}
\newtheorem*{L2}{Lemma 2}
\newtheorem*{L3}{Lemma 3}
\newtheorem*{L4}{Lemma 4}
\newtheorem*{L5}{Lemma 5}
\newtheorem*{L6}{Lemma 6}
\newtheorem*{L7}{Lemma 7}
\newtheorem*{L8}{Lemma 8}
\newtheorem*{L9}{Lemma 9}

\newtheorem*{Ex1}{Example 1}
\newtheorem*{Ex2}{Example 2}

\begin{document}
\title{Residual properties of groups defined by basic commutators}
\author{Gilbert Baumslag}
\address{CAISS and Department  of Computer Science,
City College of New York, Convent Avenue and 138th Street, New
York, NY 10031} \email{gilbert.baumslag@gmail.com}
\thanks{The research of the first author is supported by Grant CNS 111765.}

\author{Roman Mikhailov}
\address{Chebyshev Laboratory, St. Petersburg State University, 14th Line, 29b,
Saint Petersburg, 199178 Russia and St. Petersburg Department of
Steklov Mathematical Institute} \email{rmikhailov@mail.ru}
\urladdr{http://www.mi.ras.ru/\~{}romanvm/pub.html}
\thanks{The research of the second author is supported by the Chebyshev
Laboratory  (Department of Mathematics and Mechanics, St.
Petersburg State University)  under RF Government grant
11.G34.31.0026}

\begin{abstract}
In this paper we study the residual nilpotence of  groups defined by basic
commutators. We prove that the so-called Hydra
groups as well as certain of their generalizations and quotients
are, in the main, residually torsion-free  nilpotent. By way of
contrast we give an example of a group defined by two basic
commutators  which is not residually torsion-free nilpotent.
\end{abstract}
 \maketitle

\section{Introduction}

Let $\mathcal{P}$ be a property or class of groups. Then a group
$G$ is termed residually $\mathcal{P}$ if for each $g\in G, g\neq
1$, there exists a normal subgroup $N$ of $G$ such that $g \notin
N$ and $G/N \in \mathcal{P}$.  In 1935 Wilhelm Magnus \cite{WM1}
proved that free groups are residually torsion-free-nilpotent and
as a corollary that an $n<\infty$ generator group with the same
lower central quotients as a free group of rank $n$ is free. The
genesis of this paper is an earlier proof of ours, which we
include here, that the so-called Hydra groups, which are
one-relator groups defined by basic commutators and recently
introduced by Dison and Riley \cite{Hydra2},
 are residually torsion-free nilpotent. Whether all one-relator groups defined by basic commutators
are residually torsion-free nilpotent remains to be determined.
Here we prove that a number of one-relator groups defined by
 basic commutators, and some of their generalizations,
 are residually torsion-free nilpotent.   Whether residual
torsion-free nilpotence, can be used to further our understanding
of the isomorphism problem for one-relator groups seems worth
exploring further.

\subsection*{A little history} There is now a large body of work
devoted to residual properties of groups and in particular, to
residual nilpotence. Perhaps the first residual property of free
groups was obtained  by F. W. Levi \cite{FWL} in 1930 who proved
in particular that free groups are residually finite $2$-groups.
In 1935 W. Magnus \cite{WM1} proved  an even stronger theorem,
namely that free groups are residually torsion-free nilpotent.
Since K. W. Gruenberg \cite{KWG} later proved that finitely
generated torsion-free nilpotent groups are residually finite
$p$-groups for every prime $p$,  Magnus' theorem is indeed a
generalization of Levi's theorem. A. I. Malcev \cite{AIM}
subsequently extended Magnus' theorem to free products of
torsion-free nilpotent groups by proving that the free product of
residually torsion-free nilpotent groups is again residually
torsion-free nilpotent. Proceeding in a different direction, A.I.
Lichtman \cite{AIL} has studied the residual nilpotence of the
multiplicative group of a skew field generated by universal
enveloping algebras. Much work has also focussed on the residual
nilpotence of groups which are free in certain solvable varieties
of groups, including the variety of all solvable groups of a given
derived length \cite{KWG}. The proof by B. Hartley \cite{BH} that
the wreath products of torsion-free abelian groups by torsion-free
nilpotent groups are residually torsion-free nilpotent touches on
this work of Gruenberg \cite{KWG}. The techniques used by
Gruenberg as well as that of Hartley  \cite{BH} make use of basic
commutators, to be defined below, which go back to P. Hall's
fundamental work on finite groups of prime power order
\cite{PHall1} and his so-called collection process and the work of
M. Hall \cite{MHall}. A key ingredient of this work is the
introduction and use of Lie and associative rings in furthering
our understanding  of the lower central sequence of a group  that
goes back to P. Hall \cite{PHall1}, and was developed  among
others by M. Lazard \cite{MLazard} and W.  Magnus (cf.,
\cite{MKS}) and E. Witt \cite{EW} who proved that the sub-Lie
rings of a free Lie ring are again free. Here we will make much
use of the fundamental work of J.P. Labute \cite{JPL} and his
study of the Lie ring of a one-relator group, which given the
right hypothesis, turns out to be a one-relator Lie ring. A brief
survey of the residual nilpotence of a number of groups and some
additional results and references can be found in the work of the
first author \cite{Baumslagmfn}.

Much of the foregoing discussion has focussed on the purely
combinatorial group theoretic aspects of residually nilpotent
groups. They arise  naturally, not only in combinatorial group
theory but also in many geometric problems which involve knots and
links, arrangements of hyperplanes, homotopy theory, four
manifolds and many other parts of mathematics.  In particular, J.
R. Stallings \cite{JRS} proved that the lower central quotients of
a fundamental group are invariant under homology cobordism of
manifolds. His result underpins Milnor and Massey product
invariants of links, which have analogues for knots in arbitrary
3-manifolds using suitable variations of the lower central series.
The homological properties of finitely generated parafree groups,
i.e., those finitely generated residually nilpotent groups with
the same lower central quotients as free groups, play an important
role in low dimensional topology, see e.g., Cochran-Orr
\cite{CochranOrr}.  Finitely generated non-free parafree groups
exist in profusion \ \cite{GB1} and  the work of Bousfield
\cite{Bousfield} contains a large number of infinitely generated
examples.

\section{Our main results}
\subsection{The Hydra groups}
In \cite{Hydra2}, Dison and Riley introduced a family of
one-relator groups
$$
G(k,a,t)=\langle a,t\ |\ [a,\underbrace {t,\dots, t}_k]=1\rangle
\, (k\geq 1)
$$
which they termed Hydra groups. These groups are infinite cyclic
extensions of finitely generated free groups. Indeed if we put
$a_0=a$,  and $a_{i}=[a_{i-1}, t]$ for $i=1,\dots,k-1$, then the
subgroup $H$ of $G(k,a,t)$ generated by $a_0,\dots,a_{k-1}$ is a
free normal subgroup of $G(k,a,t)$ and $G(k,a,t)$ is an infinite
cyclic extension of $H$. Dison and Riley \cite{Hydra2} proved that
the subgroup of $G(k,a,t)$ generated by $a_0t,\dots,a_{k-1}t$ when
$k>1$, has extremely large distortion, by contrast with the
finitely generated subgroups of free groups which have linear
distortion. Here we will prove the following
\begin{theorem} The Hydra groups $G(k,a,t)$
are residually torsion-free nilpotent.
\end{theorem}
So the Hydra groups like free groups are residually torsion-free nilpotent.
\subsection{Generalizations of the Hydra groups}
We shall use Theorem 1 in the proof of the following two theorems, Theorem 2 and Theorem 3,
which taken together and given the right conditions, amount to
a considerable generalization of Theorem 1.

\begin{theorem} Let $X$ and
$Y$ be disjoint sets of generators and let $F$ be the free group on $Z=X \cup Y$.
Let $u$ be an element in the subgroup of $F$ generated by $X$ and $v$ an element
in the subgroup of $F$ generated by $Y$. If
$u$ and $v$ are not proper powers, then for every $k>1$
$$G(k,u,v)=\langle Z\ |\ r(u,v)=1\rangle ,$$
where $r(u,v)=[u,\underbrace{v,\dots, v}_{k}]$, is residually a finite
$p$-group for every prime $p$.
\end{theorem}
If we now ally Theorem 2 with Theorem L (see below), a theorem of J.P. Labute  \cite{JPL},  which we will describe in due course,
it is easy to deduce the following theorem, where we adopt the notation used in the formulation
of Theorem 2.
\begin{theorem}
Suppose that  $u$ and $v$ are basic commutators and that $k>1$.
Then $G(k,u,v)$ is residually torsion-free nilpotent.
\end{theorem}
Notice that we assume that $k>1$. We will consider the case where
$k=1$ separately.

\subsection{Big power groups and  residually torsion-free nilpotent groups defined by simple commutators}
Theorem 1 is at the heart of the proof of Theorem 2. Our first proof of Theorem 3 made heavy use of a rather different
result which depends on a property of free groups now termed the {\it big powers property}.  This property was used
initially to prove that certain HNN extensions of a very simple type called extensions of centralizers, are residually free in Baumslag \cite{GB0} and,  as it turns out, a little earlier by R.C. Lyndon \cite{RCL}:
\begin{theorem}
Suppose that
$$G=\langle x_1,\dots,x_n,t\ |\ [u,t]=1\rangle,$$
where $u$ is a word in $x_1,\dots,x_n$ which is not a proper power. Then $G$ is residually free.
\end{theorem}

We will still avail ourselves use of this property, which was the genesis
of our next theorem, as well as some further results which are needed to prove
Theorem 5:
\begin{theorem}
Let  $u$ be an element in the free group $F$ on $X$  and $v$ be an element in the free group $E$ on $Y$. Suppose that
$u$ and $v$ are not proper powers and  that
$$G=\langle X \cup Y\ |\ [u,v]=1\rangle.$$
Then
\begin{enumerate}
\item  $G$ is residually a finite $p$-group for every prime $p$.
\item If $u\in \gamma_j(F), u \notin\gamma_{j+1}(F)$ and  if $v\in \gamma_k(E),\, v \notin\gamma_{k+1}(E)$ and  if $u$ is not a proper power modulo $\gamma_{j+1}(F)$ and
$v$ is not a proper power modulo $\gamma_{k+1}(E)$
then $G$ is residually torsion-free nilpotent. In particular, if $u$ and $v$ are basic
commutators, then $G$ is residually torsion-free nilpotent.
\end{enumerate}
\end{theorem}
So for instance the one-relator groups
$$G(x_1,\dots,x_n)=\langle x_1,\dots,x_n\ |\ [x_1,\dots,x_n]=1\rangle$$
and
$$G=\langle a,b,c,d\ |\ [[a,b],[c,d,d]]=1\rangle$$
are residually torsion-free nilpotent. Actually more is true in
some special instances since the groups $G(x_1,\dots,x_n)$ are
even residually free.

Labute's Theorem L, already cited above, adds to what is already known about some
residual properties of certain groups, termed cyclically pinched one-relator groups,
as detailed in Subsection 2.4 below.

\subsection{Cyclically pinched one-relator groups}
A one-relator group $G$ is termed cyclically pinched if it is an amalgamated product of two free groups with
a cyclic subgroup amalgamated:
$$G=\langle x_1,\dots,x_m,y_1,\dots,y_n\ |\ u(x_1,\dots,x_m)=v(y_1,\dots,y_n)\rangle,$$
where here $u=u(x_1,\dots,x_m)$ and $v(y_1,\dots,y_n)$  are non-trivial elements  respectively in the free groups on
$X=\{x_1,\dots,x_m\}$ and $Y=\{y_1,\dots,y_n\}$. There has been considerable
attention paid to these  cyclically pinched one-relator groups since they are generalizations of the fundamental
groups of surfaces. Here we will add a little more to what is already known by invoking one of Labute's fundamental theorems in order to prove
\begin{theorem}
Let $G$ be the amalgamated product
$$G=\langle x_1,\dots,x_m,y_1,\dots,y_n\ |\ u(x_1,\dots,x_m)=v(y_1,\dots,y_n)\rangle.$$
If $u$ and $v$ are basic commutators then $G$ is residually torsion-free
nilpotent.
\end{theorem}
So for example it follows that
$$G=\langle x_1,x_2,y_1,y_2\ |\ [x_1,x_2]=[y_1,y_2,y_2]\rangle$$
is residually torsion-free nilpotent. In addition, it also follows
that if $u_1,\dots,u_q$ are basic commutators in disjoint sets of
generators coming from the set $\{x_1,\dots,x_n\}$ then the
following variation of Theorem 4 holds: the group
$$\langle x_1,\dots,x_n\ |\ [u_1,\dots, u_q]=1\rangle$$ is residually torsion-free
nilpotent. These examples should be compared with a number of
other related examples (cf., e.g., \cite{Baumslag}, \cite{ij}).

\subsection{Examples}
In the closing sections of this paper, Sections 11 and 12,  we will discuss  some
 examples defined by two basic commutators. They are all
quotients of the Hydra groups. Some of them are residually torsion-free nilpotent and some are not.
The most important and most difficult to prove example is the following, which we record as our final
theorem, in Section 11:
\begin{theorem}
The group
$$\langle a,t\ |\ [a,t,t]=[a,t,a,a,a]=1\rangle$$
is not residually torsion-free nilpotent.
\end{theorem}
We will discuss some additional examples in Section 12.

\section{The arrangement of this paper}
We will introduce some of the notation to be used here in Section
4.1. In Section 4.2 we record a little of what we will need in
order to use the Lie ring of a group and we remind the reader of
the definition of basic sequences and basic commutators in Section
4.3. Section 4.4 records some of the work  of Kim and McCarron
that we will need here,  conveniently re-stated as Theorem KM,
  and Labute's important theorem is described in Section 4.5 as Theorem L. Finally we will need
  a theorem of P. Hall which we describe in Section 4.6. The proofs of our theorems are given in
  the subsequent sections labelled by the theorems being proved ending with Section 11 where Theorem 7
  is proved, and Section 12, where some additional examples are discussed.

\section{Notation, definitions and some basic theorems of Kim and McCarron,  Labute and Philip Hall}
\subsection{Notation}
As usual, if $x,y,a_1,\dots,a_{k+1}$ are elements of a group $G$ we
set $[x,y] =x^{-1}y^{-1}xy$, $x^y=y^{-1}xy$ and define
$$[a_1,a_2,\dots, a_{k+1}]=[[a_1,\dots,a_{k}], a_{k+1}]\ (\, k>1).$$
The lower central series of $G$ is defined inductively by setting
$\gamma_1(G)=G$, $\gamma_{n+1}(G)= [G,\gamma_n(G)].$
\subsection{The Lie ring  of a group}
Each of the factor groups $L_n= \gamma_n(G)/\gamma_{n+1}(G)$ is an abelian group
and will often be written additively.  We now put
$$L(G)=\bigoplus_{n=1}^{\infty} L_n.$$
$L(G)$ can be turned into a Lie ring (over $\mathbb{Z}$) by
defining a binary operation, denoted  $[x,y]$, first on the  $L_n$
and extended by linearity to all of $L(G)$, as follows:
$$[a\gamma_{i+1}(G),b\gamma_{j+1}(G)]= (a^{-1}b^{-1}ab)\gamma_{i+j+1}(G) \, (a\in \gamma_i(G),
b\in \gamma_j(G)).$$ We term $L(G)$ the Lie ring of $G$. The
following references can be consulted by a reader interested in
the construction of such Lie rings, see for example \cite{PHall1},
\cite{MKS}.

We will need also here the definition of a basic sequence and a basic commutator.
\subsection{Basic commutators}
Let
$$X=\{x_1,\dots,x_q\}$$
be a non-empty finite set and let $G$ be the free groupoid generated by $X$.
So the elements of $G$ are simply the bracketed products of the elements of $X$
and two such products are equal only is they are identical. The number of factors  $\mid g \mid$ in
such a product $g$ is termed the length or weight of $g$. A sequence $b_1,b_2,\dots $ of elements of
$G$ is termed a basic sequence in $X$ if
\begin{enumerate}
\item every element of $X$ occurs in the sequence;
\item if $\mid b_i \mid < \mid b_j \mid$, then $i<j$;
\item if $u=vw (v,w \in G)$ is an element of $G$ of length at least 2,
$u$ occurs in the above sequence if and only if
 $v=b_i, w=b_j$ and $j<i$ and
 either $\mid b_i \mid=1$ or $b_i=b_kb_{\ell}$ and $\ell \leq j$.
\end{enumerate}
The terms in a basic sequence are called basic commutators.
The proof of the existence of such basic sequences can be found for example in
\cite{MHall}. Now if $F$ is a free group, freely generated by the set $X$
and if $\gamma_n(F)$
denotes the $n^{th}$ term of the lower central series of $F$ then Wilhelm Magnus \cite{WM1} proved
that the basic commutators of weight $n$ freely generates modulo $\gamma_{n+1}(F)$ the free abelian
group $\gamma_n(F)/\gamma_{n+1}(F)$.

\subsection{The theorems of Kim and McCarron}
We will need special cases of Theorems 3.4 and 4.2  in Kim and McCarron \cite{KM} which, for
convenience, we record here  as Theorem KM. These theorems
make use of what they call a $p$-preimage closed subgroup of a group. Here we will use a more customary notation, namely that of  a $p$-group separated subgroup. This more directly reflects
what is needed in the proof that certain amalgamated products of residually finite $p$-groups
are again residually finite $p$-groups:
\begin{definition}
A subgroup $H$ of a residually finite $p$-group $G$ is termed $p$-group separated in $G$ if for each element
$g \in G$, $g \notin H$, there exists a homomorphism $\phi$ from $G$ into a finite $p$-group such that
$\phi(g)\notin\phi( H)$.
\end{definition}
We are now in position to formulate
\begin{KM}
\begin{enumerate}
\item If $G$ is residually a finite $p$-group, if $G'$ is a copy of $G$ and if $H'$ denotes the
copy of $H$ in $G'$, then the amalgamated product $G*_{H=H'}G'$ is residually a finite
$p$-group if and only if $H$ is $p$-group separated in $G$.
\item Let $p$ be a prime, $A$ and $B$ residually finite $p$-groups and let $c\in A$, $d\in B$
be elements of infinite order. If the subgroup of $A$ generated by $c$ is $p$-group separated
in $A$ and if the subgroup of $B$ generated by $d$ is $p$-group separated in $B$ then the
amalgamated product  $A*_{c=d}B$ is residually a finite $p$-group.
\end{enumerate}
\end{KM}
\subsection{Labute's theorem}
The theorem of Labute mentioned above is then
\begin{lab} Let $F$ be a free group freely generated by $X$,  let $r$ be an element of $F$
and let $G=\langle X\ |\ r=1\rangle$. Suppose that $r \in
\gamma_n(F)$, $r \notin \gamma_{n+1}(F)$ and that
$r\gamma_{n+1}(F)$ is not a proper power (i.e., multiple) in
$L_n$. Then the Lie ring of $G$ is additively free abelian. If $G$
is residually nilpotent, then $G$ is residually torsion-free
nilpotent.
\end{lab}

\subsection{Hall's theorem}
There is yet another theorem that we will need in this paper, due to Philip Hall \cite{PHall2}:
\begin{PH}
Let $G$ be a group with a normal nilpotent  subgroup  $H$. Suppose
that $G/H$ is nilpotent  and that $H/[H,H]$ is nilpotent. Then $G$ is nilpotent.
\end{PH}

\section{The Hydra groups are residually torsion-free nilpotent}
We begin with the proof of Proposition 1.
\subsection{Proposition 1}
We will show that the proof of Theorem 1 is a consequence of  the following proposition
which is itself an easy consequence of Philip Hall's Theorem H.
\begin{P1}
Let $G$ be a group with a normal residually torsion-free nilpotent
subgroup $H$. If $G/H$ is torsion-free nilpotent and if $G/[H,H]$
is nilpotent, then $G$ is residually torsion-free nilpotent.
\end{P1}
It suffices to prove that if $g\in H, g \neq 1$, then there is
normal subgroup $K_g$ of $G$ such that $g \notin K_g$ with $G/K_g$
torsion-free nilpotent. Since $H$ is residually torsion-free
nilpotent, there is a characteristic subgroup $K_g$ of $H$  which
does not contain $g$ with $H/K_g$ torsion-free nilpotent.
Moreover, since $K_g$ is characteristic in $H$, $K_g$ is normal in
$G$.  So $G/K_g$ is an extension of a torsion-free nilpotent group
$H/K_g$ by a torsion-free nilpotent group $G/H$. Hence  $G/K_g$ is
torsion-free. Now put $\bar{G}=G/K_g$, $\bar{H}=H/K_g$. Then
$\bar{G}/[\bar{H},\bar{H}] \cong G/[H,H]$ is nilpotent by
assumption which by Hall's Theorem H, implies that $\bar{G}=G/K_g$
is nilpotent. Since $g \notin K_g$ and $G/K_g$ is torsion-free
nilpotent, this proves Proposition 1.

\subsection{The proof of Theorem 1}
We are now in a position  to prove that the Hydra groups
$$
G(k,a,t)=\langle a,t\ |\ [a,\underbrace{t,\dots, t}_k]=1\rangle \
k\geq 1
$$
are residually torsion-free nilpotent.  Since $G(1,a,t)$ is free
abelian of rank 2, it suffices to  assume that $k>1$. As we noted
above, $G=G(k,a,t)$ is an infinite cyclic extension of the free
subgroup $H$ generated by $k$ elements $a_0=a, a_i=[a_{i-1},t]$
for $i=1,\dots,k-1$.   Now free groups are residually torsion-free
nilpotent. Moreover $G/\gamma_2(H)$ is clearly nilpotent of class
$k$. So it follows by Proposition 1 that $G$ is residually
torsion-free nilpotent, as claimed.

\section{The proof of Theorem 2}
We turn now to the generalizations $G(k,u,v)$ of the Hydra groups discussed in Subsection 2.2, where
$$G(k,u,v)=\langle X,Y\ \mid r(u,v)=1\rangle,$$
$r(u,v)=[u,\underbrace{v,\dots, v}_{k}]$ and neither $u$ nor $v$ is a proper power. Our objective is to prove,
given the appropriate hypothesis, that the $G(k,u,v)$ are residually torsion-free nilpotent.  The first
step in the proof, which depends heavily on Theorem KM, which is due to Kim and McCarron \cite{KM}
and detailed in Section 4.4,  is to prove that they are residually finite $p$-groups for every
 prime $p$.  The proof  is divided
up into a number of steps which involve centralizers of elements and separation properties of
various subgroups.

\subsection{Centralizers of elements in  the Hydra groups}
We will need some more information about the Hydra groups:
\begin{L1}
The centralizer of $a$ in
$$G(k,a,t)=\langle a,t\ |\ [a,\underbrace{t,\dots, t}_{k}]=1\rangle \, (k>1)$$
 is generated by $a$.
\end{L1}
\begin{proof}
Let $H$ be the normal closure in $G=G(k,a,t)$ of $a$. Then, adopting the
notation introduced in Section 5.2 in the course of the proof of Theorem 1,
we have already noted that $H$ is free on the $a_j$ and that
$G$ is the semi-direct product of $H$ and the infinite cyclic group on $t$. In addition,
 $t$ acts on $H$ as follows:
$$t^{-1}a_0t=a_0a_1,\ \dots,\ t^{-1}a_{k-2}t=a_{k-2}a_{k-1},\ t^{-1}a_{k-1}t=a_{k-1}.$$
Note that $a=a_0$.

Suppose that $g \in G(k,t,a)$ and that $[g,a]=1$. Since $G(k,t,a)=H \rtimes
\langle t \rangle$, $g=ht^n\,  (h\in H)$. We can assume that
$n\geq 0$. Then
 $$a_0=g^{-1}a_0g=t^{-n}h^{-1}a_0ht^n.$$
 Now working modulo $\gamma_2(H)$, we  find $a_0=t^{-n}a_0t^n$. It follows, again
 working modulo $\gamma_2(H)$, that
 $$t^{-n}a_0t^n=a_0a_1^n\dots, $$
 where the terms following $a_1^n$ are words in $a_2,\dots,a_{k-1}$.
 Hence these words are trivial in $H$ and $n=0$. It follows that we have proved that
 $$a_0=h^{-1}a_0h.$$
 But $H$ is free and therefore the centralizer of $a_0$ in $H$ is a power of $a_0$. This completes the
 proof.
 \end{proof}

 Now denote the free group on $Y$ by $F(Y)$. Then one of the consequences of Lemma 1 is
\begin{L2}
Let $k\geq 2$,  let $v\in F(Y)$ be an element which is
not a proper power in $F(Y)$ and let
$$
J=\langle Y,a,t\mid [a,\underbrace{t,\dots, t}_k]=1,t=v\rangle .$$
Then the  centralizer of $a$ in $J$ is the cyclic subgroup
$\langle a\rangle$.
\end{L2}

\begin{proof}
Observe that  $J$ is an amalgamated product of $G(k,t,a)$ and $F(Y)$:
$$J=G(k,t,a)*_{t=v} F(Y).$$
Suppose that $g$ lies in the centralizer of $a$ in $J$.
Then $g$ can be written in the form
$$
g=e_1\dots e_n,$$ where each $e_i$ is from one of the factors
$F(Y)$ or $G(k,a,t)$, and successive $e_i,e_{i+1}$ come from different
factors. Since $ag=ga$ and $a$ does not belong to the subgroup
generated by $t$ it is not hard to see that $n=1$. In this case,
we have $ae_1=e_1a$. It follows that $e_1 \in G(k,t,a)$ in which case
Lemma 1 applies which means that $e_1$ is a power of
$a$ as claimed.
\end{proof}

\subsection{$p$-group separated subgroups}
We will need now to make use of  the work of  Kim and McCarron on residually finite $p$-groups.
Our objective is to use the fact that   $G(k,u,v)$ is an amalgamated product in which the amalgamated
subgroups satisfy a separation property. This is accomplished by the following lemma.

\begin{L3}
The subgroup $\langle t \rangle$ is $p$-group separated in $G(k,a,t)$
for all primes $p$.
\end{L3}

\begin{proof}
Finitely generated, residually torsion-free nilpotent groups are
residually finite $p$-groups for every choice of the prime $p$
\cite{KWG}. So $G(k,a,t)$ is residually a finite $p$-group by Theorem 1.
Consequently,  it follows from Kim and McCarron's Theorem KM as described in
Subsection 4.4, that
$\langle
a\rangle$ is $p$-group separated in $G(k,a,t)$ if and only if the group
$$
\tilde G(k,a,t):=\langle a,t,c\ |\ [a,\underbrace{t,\dots,
t}_{k}]=[c,\underbrace{t,\dots, t}_{k}]=1\rangle
$$
is residually a finite $p$-group. Now $\tilde G(k,a,t)$ is the middle of a short
exact sequence
$$ 1\to H_{2k}\to \tilde G(k,a,t)\to
\mathbb \langle t\rangle\to 1
$$
where  $H_{2k}$ is the free subgroup generated by
$\{a,a^t,\dots,a^{t^{k-1}},c,c^t,\dots, c^{t^{k-1}}\},$ and the
quotient of $\tilde G(k,a,t)$ by $H_{2k}$ is  cyclic and generated
by the image of $t$. It follows along the same lines as in  the
proof of Theorem 1 that the group $\tilde G(k,a,t)$, with a slight
abuse of notation, is residually the semi-direct product of the groups
$(H_{2k}/\gamma_i(H_{2k})\rtimes \langle t\rangle$ $(i\ge
1)$. Moreover, the groups $H_{2k}/\gamma_i(H_{2k})\rtimes \langle
t\rangle$ are finitely generated, torsion-free nilpotent for all
$i\geq 1$. Hence $\tilde G(k,a,t)$ is residually a finite
$p$-group and therefore $\langle t\rangle$ is $p$-subgroup
separated in $G(k,a,t)$, as claimed.
\end{proof}

Next we have
\begin{L4}
Let
$$
J=\langle Y,a\ |\ [a,\underbrace{v,\dots, v}_{k}]=1\rangle,
$$
where $v$ is not a proper power in F(Y). Then $J$ is residually
a finite $p$-group for every prime $p$.
\end{L4}

\begin{proof}
Observe that $J$ is an amalgamated product:
$$
J=F(Y)*_{v=t}G(k,a,t).
$$
In order to prove that  $J$ is residually a finite $p$-group we
again have to appeal to Theorem KM. We have already
proved that
 $\langle t\rangle$ is $p$-subgroup separated in
 $G(k,a,t)$. Now by Theorem 4, since $v$ is not a
proper power, $L=\langle Y,t\ |\ [t,v]=1\rangle$ is residually
free. So the subgroup $M$ of $L$ generated by $Y$ and $t^{-1}Yt$
is residually a finite $p$-group. But $M$ is an amalgamated
product of two copies of $gp(Y)$ amalgamating $gp(v)$. Hence
$gp(v)$ is $p$-group separated in $gp(Y)$. It follows that $J$ is
residually a finite $p$-group by Theorem KM.
\end{proof}

The last step before we come to the proof of Theorem 2 is
\begin{L5}
The subgroup $\langle a\rangle$ is $p$-subgroup separated  in $J$.
\end{L5}
\begin{proof}
Suppose that $c\in J$ is such that for every homomorphism $\phi$ from
$J$ into a finite $p$-group, $\phi(c) \in \phi(\langle a \rangle)$. It
follows that $[c,a]$ lies in every normal subgroup of $J$ of index
a power of $p$. Since $J$ is residually a finite $p$-group, it
follows that $[c,a]=1$. So if $k>1$, $c$ is a power  of $a$ since the
centralizer of $a$ in $J$ is generated by $a$.
\end{proof}

We are now in position to complete the proof of Theorem 2.
Present $G_k(u,v)$ as an amalgamated product
$$
G(k,u,v)=J*_{a=u}F(X).
$$
Since $J$ is residually a finite $p$-group, $\langle a\rangle$
is $p$-subgroup separated in $J$, and $\langle u\rangle$ is
$p$-subgroup separated in $F(X)$, the group $G(k,u,v)$ is
residually a finite $p$-group by Theorem KM.

\section{The proof of Theorem 3}
The proof of Theorem 3 follows almost immediately from Labute's
Theorem L. In order to explain why, suppose now that $F$ is the
free group on $Z=X\cup Y$. Because $u$ and $v$ are basic
commutators, say of weights $m$ and $n$, respectively,
 in disjoint sets of generators $r(u,v)$ is a basic commutator of weight $m+kn$. Consequently
 $r(u,v) \in \gamma_{m+kn}(F)$ and $r(u,v) \notin \gamma_{m+kn+1}(F)$. Moreover
$r(u,v)  \gamma_{m+kn+1}(F)$ is an element in a basis for $ \gamma_{m+kn}(F)/\gamma_{m+kn+1}(F)$
and hence is not a proper power. Since $G(k,u,v)$ is residually a finite $p$-group, by Theorem 2,
it follows by Theorem L, that $G(k,u,v)$ is residually torsion-free nilpotent.

\section{The proof of Theorem 4}
Suppose that $u$ is an element of a free group $F$ that is not a
proper power and that $G$ is the one-relator group on $X \cup
\{t\}$ defined by the single relation $[u,t]=1$. Then $G$ is
residually free by \cite{GB1}. Now $G(x_1,\dots,x_n)$ is simply
obtained from the free group on $x_1,\dots,x_{n-1}$ by adding an
additional generator $x_n$ and the relation
$[[x_1,\dots,x_{n-1}],x_n]=1$. By a theorem of Magnus, Karrass and
Solitar \cite{MKS} a non-trivial commutator in a free group is not
a proper power. Hence Theorem 4 is an immediate consequence of
\cite{GB0}.

\section{The proof of Theorem 5}
It turns out that the main step in the proof of Theorem 5 is the following lemma.
\begin{L6}
Let $G$ be the one-relator group
$$G=\langle X\cup \{t\}\ |\ [u,t]=1\rangle,$$
where $u$ is an element in the free group $E$ on $X$ which is not a proper power. Then the
following hold:
\begin{enumerate}
\item $G$ is residually a finite $p$-group for each prime $p$;
\item the subgroup of $G$ generated by $t$ is finitely $p$-group separable in $G$.
\end{enumerate}
\end{L6}
That $G$ is residually a finite $p$-group  follows, as previously noted, because $G$ is even residually
free. However in order to prove (2), we need some additional information.
 To this end, notice  that if $A$ is the free abelian group on $s$ and $t$,
if $H$ is the subgroup of $E$ generated by $u$, $K$ is the subgroup of $A$ generated by $s$,
then $G$ can be viewed as an amalgamated product of $E$ and $A$ with $H$ amalgamated with
$K$ according to the the isomorphism $\phi$ mapping $u$ to $s$:
$$G=\langle E * A\ | H=_{\phi} K\rangle.$$
Now free groups and free abelian groups
are residually finite $p$-groups. Moreover, $H$ is finitely $p$-group separable
in $E$ since $u$ generates its centralizer in $E$ and $K$ is
clearly finitely $p$-group separable in $A$. So it follows from Theorem KM that $G$ is residually a finite
$p$-group.  This again  proves (1).

We are left with the proof of (2), that is if $g\in G$ and $g
\notin T=\langle t\rangle$, then there is a homomorphism $\phi_g$
of $G$ into a finite $p$-group such that $\phi_g(g) \notin \phi_g(T)$.
The proof will be divided up into a number of cases.

\begin{enumerate}
\item  $g\in H, g \neq 1$. Of course, $g \notin T$. Since $E$ is
residually a finite $p$-group, there exists a normal subgroup $I$
of $E$ of index a power of $p$ such that $g \notin I$. Define
$\phi_g:G \longrightarrow E/I$ which maps $E$ onto $E/I$ and $A$
onto $E/I$ by sending $s$ onto $uI$ and $t$ to the identity. Then
$\phi_g(g) \notin \phi_g(T)$, as required. \item $g\in E$, $g
\notin H$. Then $[g,u] \neq 1$. Choose a normal subgroup $I_1$ of
$E$ of index a power of $p$ such that $[g,u] \notin I_1$. Then
$gI_1 \notin TI_1$ for otherwise $[g,u]\in I_1$. Then $uI_1$ is a
non-trivial element of $E/I_1$ of order a power of $p$.  $\phi_g:G
\longrightarrow E/I_1$ is defined first on $E$ and then on $A$. We
take it to be the canonical homomorphism of $E$ onto $E/I_1$. Next
we define $\phi_g$ to be  the homomorphism of  $A$ to $G/I_1$
which maps $s$ to $uI_1$ and $t$ to the identity. Let $I$ be the
kernel of $\phi_g$. Then $\phi_g(g) \neq 1$ and $\phi_g(T)
=\{1\}$. So $\phi_g(g) \notin \phi_g(T)$ as needed. \item $g \in
A, g \notin T$. Then $g=s^kt^\ell$ where $k\neq 0$. Since $E$ is
residually a finite $p$-group we can choose a homomorphism
 of $E$ into a finite $p$-group so that the image of $u^k$ has arbitrarily large order
a power of $p$.  In addition, there exists a homomorphism of $A$
into a finite group so that the image of $s$ has arbitrarily large
finite order and the image of $t$ is 1. It follows that there
exists a homomorphism $\phi_g$ of $G$ into a finite $p$-group $P$
which maps $T$ to $\{1\}$ and $g$ to an element outside the image
of $T$. \item $g=f_1a_1\dots f_na_n$ where the $f_j \in E, f_j
\notin H$, $a_j \in A, a_j \notin K$. Notice that if $g \notin A$,
then
$$tg= tf_1a_1\dots f_na_n\neq f_1a_1\dots f_na_nt,$$
because $f_1 \notin A$. Now $G$ is residually a finite $p$-group and $[g,t]\neq 1$. Hence there is
a homomorphism $\phi_g$ of $G$ into a finite $p$-group which maps $[g,t]$ to a non-trivial
element. So $\phi_g(g) \notin \phi_g(T)$ since this implies that  $\phi_g([g,t])=1$.
\end{enumerate}
Since $G(u,v)$ is an amalgamated product of two residually finite
$p$-groups  where the amalgamated subgroups are cyclic and
finitely $p$-group separable, $G(u,v)$ is residually a finite
$p$-group. Now suppose that $u \in \gamma_j(F), u \notin
\gamma_{j+1}(F)$, $u\gamma_{j+1}(F)$  is not a proper power and
that $v \in \gamma_k(F), u \notin \gamma_{k+1}(F)$,
$v\gamma_{k+1}(F)$ is not a proper power, then it is not hard to
prove that $[u,v] \in \gamma_{j+k}(F), [u,v] \notin
\gamma_{j+k+1}(F)$ and that  $[u,v] \gamma_{j+k+1}(F)$ is not a
proper power. So Theorem L applies and therefore $G(u,v)$ is
residually torsion-free nilpotent.

\section{The proof of Theorem 6}
Our objective now is prove Theorem 6, namely that if
$$G=\langle x_1,\dots,x_m,y_1,\dots,y_n\ |\ u(x_1,\dots,x_m)=v(y_1,\dots,y_n)\rangle$$
and if $u$ and $v$ are basic commutators, then $G$ is  residually torison-free  nilpotent.
Now Magnus, Karrass and Solitar \cite{MKS4} have proved that in a
free group, a non-trivial commutator
is not a proper power. Consequently $G$ is residually a
 finite $p$-group \cite{Baumslag}. In addition, $uv^{-1}$ is then a
product of a basic commutator of weight $k$, say, and the inverse of a
second basic commutator of weight $\ell$,
say, in a disjoint set of generators. So if $q$ is the minimum of $k$ and $\ell$, then $uv^{-1}\in \gamma_q(F),
uv^{-1} \notin \gamma_{q+1}(F)$, where $F$ is the free group on $x_1,\dots,x_m,y_1,\dots,y_n$ and $uv^{-1}\gamma_{q+1}(F)$
is not a proper power in $\gamma_q(F)/\gamma_{q+1}(F)$. So by Labute's Theorem L,  $G$ is residually torsion-free nilpotent.

\section{The proof of Theorem 7}
In this and the subsequent section we shall use Proposition 1 to define several examples of quotients of the Hydra groups
which are residually torsion-free nilpotent. We begin first with the proof of Theorem 7, namely that  the group
$$D=\langle a,t\ |\ [a,t,t]=[a,t,a,a,a]=1\rangle$$
is not residually  a finite p-group if $p\neq 2$. As noted previously,  Gruenberg \cite{KWG} proved that
finitely generated torsion-free nilpotent groups are residually finite $p$-groups for all primes $p$.
So $D$ is not residually torsion-free nilpotent.
The proof is complicated and will be carried out by means of a series of
lemmas.

Our proof uses the well-known Hall-Witt identity several times.
Recall it for the readers convenience. For elements $A,B,C$ of a
group, the following identity is satisfied:
\begin{equation}\label{hw}
[A,B^{-1},C]^B[B,C^{-1},A]^C[C,A^{-1},B]^A=1.
\end{equation}

In the rest of this section, we will use the notation
$w=[a,t,a,a,t,a]$ . We will show that $w\neq 1$ is a
generalized 2-torsion element of $G$. That is, for every $n\geq
1$, the order of $w$ is a power of 2 modulo $\gamma_n(G)$.

\begin{L7}
For every $c\in [G,G]$ and $n\geq 2$,
$$
[a,t,a,a,t,c]\in [\langle w\rangle^G,G]\gamma_n(G),
$$
where $\langle w\rangle^G$ is the normal closure of $w$ in $G$.
\end{L7}
\begin{proof}
We will prove the statement for $c=[a,t]$. The general case
clearly will follow, since the general element of $[G,G]$ can be
presented as a product
$$
\prod_j [a,t]^{\pm a^{l_j}q_j},\ l_j\in \mathbb Z,\ q_j\in [G,G].
$$

Denote $v:=[a,t,a,a,t,[a,t]]$. The Hall-Witt identity (see
(\ref{hw}) with $A=[a,t,a,a], B=t^{-1},C=[a,t]$) implies that
$$
[a,t,a,a,t,[a,t]]^{t^{-1}}[t^{-1},[t,a],[a,t,a,a]]^{[a,t]}[a,t,[a,t,a,a]^{-1},t^{-1}]^{[a,t,a,a]}=1.
$$
The relation $[a,t,t]=1$ implies that $[t^{-1},[t,a]]=1$ and
therefore,
\begin{equation}\label{t4}
[a,t,a,a,t,[a,t]]=[a,t,[a,t,a,a]^{-1},t^{-1}]^{-[a,t,a,a]t}.
\end{equation}
Applying the Hall-Witt identity one more time (see (1) with $A=a,
B=t^{-1}, C=[a,t,a,a]^{-1}$), we get
\begin{equation}\label{t5}
[a,t,[a,t,a,a]^{-1}]^{t^{-1}}[t^{-1},[a,t,a,a],a]^{[a,t,a,a]^{-1}}[[a,t,a,a]^{-1},a^{-1},t^{-1}]^a=1.
\end{equation} The relation $[a,t,a,a,a]=1$
together with (2) implies that
\begin{multline*}
[a,t,[a,t,a,a]^{-1}]=[t^{-1},[a,t,a,a],a]^{-[a,t,a,a]^{-1}t}=\\
[a,t,a,a,t,a^t]^{-t^{-1}[a,t,a,a]^{-1}t}=[a,t,a,a,t,a[a,t]]^{-t^{-1}[a,t,a,a]^{-1}t}\in
\langle w\rangle^G\langle v\rangle^G.
\end{multline*}
Now the identity (\ref{t4}) implies that
$$
v\in [\langle w\rangle^G,G][\langle v\rangle^G,G]
$$
and the needed statement follows by induction on $n$.
\end{proof}

\begin{L8}  For every $n\geq 2$,
$[a,t,a^{-1},[a,t,a]]\in [\langle w\rangle^G,G]\gamma_n(G).$
\end{L8}
\begin{proof}
First observe that
\begin{multline}\label{longrel1} [a,t,a^{-1},[a,t,a]]=[[a,[a,t]]^{a^{-1}},[a,t,a]]=\\
[a,[a,t],a^{-1},[a,t,a]]=[[a,t,a,a]^{a^{-1}},[a,t,a]]^{[a,[a,t]]}=[[a,t,a,a],[a,t,a]]^{[a,[a,t]]}.
\end{multline}
Denote $E=[a,t,a,a]$. The relation $[E,a]=1$ implies that $[E^{\pm
1}, a^{\pm 1}]=1$. The Hall-Witt identity (see (\ref{hw}) with
$A=[a,t], B=a^{-1}, C=E$) implies that
$$
[a,t,a,E]^{a^{-1}}[a^{-1},E^{-1},[a,t]]^E[E,[t,a],a^{-1}]^{[a,t]}=1.
$$
Hence
$$
[E,[a,t,a]]=[E,[t,a],a^{-1}]^{[a,t]a}.
$$
Applying the Hall-Witt identity one more time, we get
$$
[E,[a,t,a]]=[[E,t^{-1},a^{-1}]^{ta},a^{-1}]^{[a,t]a}=[[t,E]^{t^{-1}},a^{-1}]^{ta},a^{-1}].
$$
The needed statement follows from Lemma 7 and the identity
(\ref{longrel1}).
\end{proof}
\begin{L9}
For every $n$,
$$
[t,a,[a,t,a]^{-1},a]\equiv w^{-1}\mod [\langle
w\rangle^G,G]\gamma_n(G).
$$
\end{L9}

\begin{proof} The Hall-Witt identity (see \ref{hw} with $A=t,B=a^{-1},C=[a,[a,t]]$) implies that
\begin{equation}\label{pa1}
[t,a,[a,[a,t]]]^{a^{-1}}[a^{-1},[a,t,a],t]^{[a,t,a]^{-1}}[a,[a,t],t^{-1},a^{-1}]^t=1.
\end{equation}
Since $[a,t,a,a,a]=1,$ one has
$$
[a^{-1},[a,t,a],t]=[a,t,a,a,t].
$$
By Lemma 7, for every $n\geq 2$,
$$
[[a,t,a,a,t]^{[a,t,a]^{-1}},a]\equiv w\mod [\langle
w\rangle^G,G]\gamma_n(G).
$$
Now the identity (\ref{pa1}) implies that the needed statement of
Lemma is equivalent to the statement that, for every $n\geq 2$,
\begin{equation}\label{wqw1}
[[a,[a,t],t^{-1},a^{-1}]^t,a]\in [\langle
w\rangle^G,G]\gamma_n(G).
\end{equation}
The Hall-Witt identity (see (\ref{hw}) with
$A=a,B=[t,a],C=t^{-1}$) and the relation $[a,t,t]=1$ imply that
$$
[a,[a,t],t^{-1}]^{[t,a]}[t,a,a^{-1},[t,a]]^{t^{-1}a}=1.
$$
We can rewrite the last identity as
$$
[a,t,a,t]^d[[a,t],[a,t,a]]^e=1,
$$
where $d=t^{-1}[a,[a,t]][t,a],\ e=[a,t,a]a[t,a]^2t^{-1}a$. The
Hall-Witt identity (see (\ref{hw}) with $A=a,B=t^{-1}, C=[a,t,a]$)
implies that
$$
[[a,t],[a,t,a]]^{t^{-1}}[t^{-1},[a,t,a]^{-1},a]^{[a,t,a]}[a,t,a,a^{-1},t^{-1}]^a=1.
$$
Therefore, there exist elements $c_1,c_2,c_3\in G$, such that
\begin{equation}\label{mz1}
[a,t,a,t]=[a,t,a,t,a^{c_1}]^{-c_2}[a,t,a,a,t]^{c_3}.
\end{equation}
By Lemma 7, for every $n\geq 2$ and $g_1,g_2\in \langle
a\rangle^G$,
$$
[[a,t,a,a,t],g_1,g_2]\in[\langle w\rangle^G,G]\gamma_n(G).
$$
Now the identity (\ref{mz1}) implies that, for every $n\geq 2$ and
$g_1,g_2\in \langle a\rangle^G$,
$$
[a,t,a,t,g_1,g_2]\in[\langle w\rangle^G,G]\gamma_n(G)
$$
and the needed statement (\ref{wqw1}) follows.
\end{proof}

We are now in a position to complete the proof of Theorem 7. We claim
that if $w=[a,t,a,a,t,a]$, then
$$w\notin
\gamma_7(G),\ w^2\in \gamma_7(G).$$
This can be proved directly or by appealing to GAP, as the following GAP fragment
shows:
\begin{align*} & \text{gap}>\text{F:=FreeGroup(2);;}\\
& \text{a:=F.1;; t:=F.2;;}\\
& \text{gap}>\text{G:=F/[ LeftNormedComm([a,t,t]),
LeftNormedComm([a,t,a,a,a])];;}\\
& \text{gap}>\text{phi:=NqEpimorphismNilpotentQuotient(G,6);;}\\
& \text{gap}>\text{aa:=Image(phi,G.1);;}\\
& \text{tt:=Image(phi,G.2);;}\\
& \text{gap}>\text{xx:=LeftNormedComm([aa,tt,aa,aa,tt,aa]);;}\\
& \text{gap}>\text{Order(xx);}\\ & 2
\end{align*}
Now we will show that $w$ is a generalized 2-torsion element. That
is, for every $n\geq 1$, the order of $w$ is a power of 2 modulo
$\gamma_n(G)$.

The Hall-Witt identity (see (\ref{hw}) with $A=[a,t,a,a],
B=t^{-1},C=a$) together with the relation $[a,t,a,a,a]=1$ implies
that
$$
[a,t,a,a,t,a]^{t^{-1}}[t^{-1},a^{-1},[a,t,a,a]]^a=1.
$$
Hence,
\begin{equation}\label{da}
[a,t,a,a,t,a]^{t^{-1}}[[a,t^{-1}],[a,t,a,a]]=1.
\end{equation}
The relation $[a,t,t]=1$ implies that
$$
[a,t^{-1}]=[t,a]^{t^{-1}}=[t,a].
$$
The identity (\ref{da}) can be rewritten as
\begin{equation}\label{da11}
[a,t,a,a,t,a]^{t^{-1}}[[t,a],[a,t,a,a]]=1.
\end{equation}
It follows from the Hall-Witt identity (see \ref{hw} with
$A=[a,t,a], B=a^{-1}, C=[t,a]$) that
\begin{equation}\label{3term}
[a,t,a,a,[t,a]]^{a^{-1}}[a^{-1},[a,t],[a,t,a]]^{[t,a]}[t,a,[a,t,a]^{-1},a^{-1}]^{[a,t,a]}=1.
\end{equation}

The second term of the relation (\ref{3term}) lies in $[\langle
w\rangle^G,G]\gamma_n(G)$ for every $n\geq 2$, by Lemma 8. The
third term of the relation (\ref{3term}) is equivalent to $w$
modulo $[\langle w\rangle^G,G]\gamma_n(G)$ for every $n\geq 2$, by
Lemma 9. Now the relations (\ref{da11}) and (\ref{3term}) imply
that, for every $n$,
$$
w^2 \equiv [\langle w\rangle^G,G]\gamma_n(G).
$$
$\Box$

\section{Two more examples}
\begin{Ex1}
Let
$$G_k=\langle a,t\ | [a,\underbrace{t,\dots,t}_k]=1,
[a,\underbrace{t,\dots,t}_{k-2},\underbrace{e,\dots,e}_{\ell}]=1\rangle
,$$ where $e=[a,\underbrace{t,\dots,t}_{k-1}]\rangle$. Then $G_k$
is residually torsion-free nilpotent for every $\ell\geq 1$.
\end{Ex1}

It is worth noting that the simplest non-nilpotent group of such
type is the following:
$$
\langle a,t\ |\ [[a,t],[a,t,t]]=1, [a,t,t,t]=1\rangle.
$$
To show that the above groups satisfy the hypothesis  of
Proposition 1, consider the following generators of the normal
closure of $a$:
\[
c_1=a,\ c_2=[a,t],\ \dots, c_j=[a,\underbrace{t,\dots, t}_{j-1}],\
j=1,\dots, k.
\]
The action of $\langle t\rangle$ on these generators is given as
already discussed in the proof of Theorem 1. Writing the relator
$[[a,\underbrace{t,\dots, t}_{k-2}],\underbrace{e,\dots,
e}_{\ell}]$ in terms of the generators $c_1,\dots, c_k$, we find
that it is $[c_{k-1}, \underbrace{c_{k},\dots, c_{k}}_{\ell}]$.
Recall that a free product of residually torsion-free nilpotent
groups is residually torsion-free nilpotent \cite{AIM}. The group
$H=\langle a\rangle^{G_k}$ (using the notation in Proposition 1)
is the free product
$$
F(c_1,\dots, c_{k-2})*\langle c_{k-1},c_k\ |\ [c_{k-1},
\underbrace{c_{k},\dots, c_{k}}_{l}]=1\rangle
$$
which is residually torsion-free nilpotent by Theorem 1. The group
$G_k$ is residually torsion-free nilpotent by Proposition 1.

Next we have

\begin{Ex2}  For $k,s\geq 1$, the
group
\begin{equation}\label{newe}
\langle a,t\ |\ [[a,\underbrace{t,\dots,
t}_{s}],[a,\underbrace{t,\dots, t}_{k-1}]]=1,
[a,\underbrace{t,\dots, t}_{k}]=1\rangle.
\end{equation}
is residually torsion-free nilpotent.
\end{Ex2}

Again, denoting the images of the $c_i$-s as before simply as $c_i$,
 we see that the subgroup $H$ which is the normal closure of the element $a$ can be presented in the form
$$
\langle c_1,\dots, c_k\ |\ [c_i,c_k]=1,\ i=s+1,\dots, k-1\rangle
$$
which is isomorphic to the group $F(c_1,\dots,
c_{s})*(F(c_{s+1},\dots, c_{k-1})\times\langle c_k\rangle )$,
which is clearly residually torsion-free nilpotent. Conditions of
the Proposition 1 are satisfied, hence the group (\ref{newe}) is
residually torsion-free nilpotent. A simple example of a group of
this kind is
$$
\langle a,t\ |\ [[a,t],[a,t,t,t]]=1, [a,t,t,t,t]=1\rangle.
$$

\vspace{.3cm}\noindent{\bf Remark.} Observe that, for $k\geq 1$,
the groups
$$
\langle a,t\ |\ [a,\underbrace{t,\dots, t}_{k-1},a]=1,
[a,\underbrace{t,\dots, t}_{k}]=1\rangle
$$
are residually nilpotent by the following result from
\cite{Mikhailov}: any central extension of a one-relator
residually nilpotent group is residually nilpotent.

As we noted at the outset, we have been unable to determine
whether a one-relator group defined by a basic commutator is
residually torsion-free nilpotent, or residually a finite p-group
or even residually finite. The best that we have managed to find
is an example of a group defined by two relations, which are basic
commutators, which is not even residually a finite p-group.

\vspace{.3cm}\noindent{\bf Remark.} For a free group on two
generators, the seventh term of the lower central series is the
normal closure of all basic commutators of weight 7,8,9,10 (see
\cite{JGS}). It follows from the proof of Theorem 7 that the group
$$
\langle a,t\ |\ [a,t,t]=[a,t,a,a,a]=1,\ \text{all basic
commutators of weight 7,8,9,10}\rangle
$$
has 2-torsion, namely
$$
[a,t,a,a,t,a]^2=1,\ [a,t,a,a,t,a]\neq 1.
$$

\vspace{.3cm} \noindent{\bf Acknowledgement.}  The authors thank
the referee for comments and suggestions. The second author thanks
A. Talambutsa for helping with GAP.

\end{document}